\documentclass[11pt]{amsart}
 \usepackage{amsopn}
 \usepackage{amsmath,amsthm,amssymb}

\usepackage{amssymb}
\usepackage{graphicx}
\usepackage{amsmath}
\usepackage[latin1]{inputenc}
\usepackage{amsfonts}

 \theoremstyle{plain}

  \theoremstyle{definition}

 \theoremstyle{remark}

\newcommand{\rr}{\mathbb{R}}

\newcommand{\nul}{\mathcal{N}}

\begin{document}
 
  \title{On completeness of integral manifolds of nullity distributions}
\author[C. Olmos and F. Vittone]{Carlos Olmos and  Francisco Vittone}

\begin{abstract}
We give a conceptual proof of the fact that if $M$ is a complete submanifold of a space form, then the maximal integral manifolds of the nullity distribution of its second fundamental form through points of minimal index of nullity are complete.
\end{abstract}

\maketitle
\section{Introduction} Let  $M$ be a submanifold of a space form and let $\nul$ be the nullity distribution of its second fundamental form.  The index of nullity of $M$ at $p$ is the dimension of $\nul_{p}$. It is well known, from Codazzi equation, that $\nul$ is an autoparallel distribution restricted to the open and dense subset of  $M$ where the index of nullity is locally constant. 

If $M$ is complete and one restricts to the open subset $U$ of points of $M$ where the index of nullity is minimal, then the integral manifolds of $\nul$ through points of $U$ are also complete, from a result of Ferus \cite{ferus}.  We will give a conceptual proof of this result, as a corollary of a general theorem, whose proof involves very simple geometric ideas.

\section{Main results} 

\textbf{Lemma.}
 \textit{Let $M$ be a Riemannian manifold and let $f:M\to N$ be a differentiable function of constant rank such that $f(M)$ is an embedded submanifold of $N$ (this can always be assumed locally). Assume that the distribution $\ker(df)$ is autoparallel and let $\Sigma$ be an integral manifold of $\ker(df)$. Let $\gamma(t)$ be a geodesic in $\Sigma$, $v\in T_{f(\gamma(0))} f(M)$ and let $J(t)$ be the horizontal lift of $v$ along $\gamma$, i.e., $J(t)\in \ker (df_{\gamma(t)})^{\bot}$ and $df J(t)=v$. Then $J(t)$ is a Jacobi vector field along $\gamma$.}
  \label{lema1}
  
  \begin{proof}
Let $c(s)$ be a (short) curve in $f(M)$ such that $c'(0)=v$. 
The horizontal lift of $c(s)$ through points of $\Sigma$ gives rise to a perpendicular variation $\widetilde{c}_{q}(s)$ by totally geodesic submanifolds, which must be by isometries. Therefore $\widetilde{c}_{\gamma(t)}(s)$ is a variation by geodesics whose associated variation field is $J(t)$. 
  \end{proof}
  
\textbf{Theorem.}
\textit{Let $M$ be a complete Riemannian manifold, $f:M\to N$ be a differentiable function and let $U$ be the open subset of $M$ where the rank of $f$ is maximal. Assume that $\ker(df)_{|U}$ is autoparallel. Then its integral manifolds are complete. }

\begin{proof}
Let $\Sigma$ be a totally geodesic integral manifold of $\ker(df)$ through a point $p\in U$ and let $\gamma:[0,b)\to \Sigma$ be a maximal geodesic in $\Sigma$. 

Observe that $f$ has maximal rank in a neighborhood of each point of $\Sigma$. From the local form of maps of constant rank it is not difficult to see that given $t_{1},\; t_{2}\in[0,b)$ there are open neighborhoods $V_{1}$ and $V_{2}$ of $\gamma(t_{1})$, $\gamma(t_{2})$ such that $f(V_{1})$ and $f(V_{2})$ are embedded submanifolds of $N$ and  $f(V_{1})\cap f(V_{2})$ contains an open neighborhood of $f(\gamma(t_{1}))=f(\gamma(t_{2}))$ in both $f(V_{1})$ and $f(V_{2})$. In particular, $T_{f(\gamma(t_{1}))}f(V_{1})=T_{f(\gamma(t_{2}))}f(V_{2})=:\mathbb{V}$.

Let $v\in \mathbb{V}$ and apply the previous lemma  to define a Jacobi field $J$ along $\gamma$ that projects down to $v$. Since $M$ is complete $\gamma(b)$ and $J(b)$ are well defined and $J(b)$, by the continuity of $df$, also projects down to $v$. Then $df_{\gamma(b)}(T_{\gamma(b)}M)$ contains $\mathbb{V}$. So, $\text{rank}(df_{\gamma(b)})=\text{rank}(df_{\gamma(0)})$ and therefore $\gamma(b)\in \Sigma$.  
\end{proof}

\medskip

If $M^{n}$ is a submanifold of the Euclidean space $\rr^{n+k}$, the Gauss map of $M$ is the map $G:M\to G_{k}(\rr^{n+k})$ definded by $p\mapsto \nu_{p}M$, where $\nu_{p}M$ denotes the normal space of $M$ at $p$. If $M^{n}$ is a submanifold of the sphere $\mathbb{S}^{n+k}\subset \rr^{n+k+1}$, then the Gauss map of $M$ is defined to be the map $G:M\to G_{k}(\rr^{n+k+1})$ that sends each point to its normal space in the sphere, regarded as a subspace of $\rr^{n+k+1}$ (see the remark bellow). A similar construction can be made for a submanifold $M^{n}$ of the hyperbolic space $H^{n+k}$, regarded as a submanifold of the Lorentz space $\rr^{n+k,1}$.

It is well known that in the three cases, the nullity distribution of $M$ coincides with the kernel of its Gauss map. Therefore we get,

\medskip
\textbf{Corollary.}
\textit{Let $M$ be a complete submanifold of a space form. Then any maximal integral manifold of the nullity distribution of $M$ through a point of minimal index of nullity is complete. }
\medskip

\textbf{Remark.}
Let $M$ be a submanifold of a space form and let $\nu_{1}$ be a parallel sub-bundle of the normal bundle. One can regard to the nullity of the second fundamental form projected to this sub-bundle. This is equivalent to regard the common kernell of all shape operators of vectors in this sub-bundle. This generalized nullity space coincides with the kernell of the generalized Gauss map $p\mapsto \nu_{1}(p)$. So, as in the corollary, if $M$ is complete one has completeness of the integral manifolds where the kernel has minimal dimension. 

\bibliographystyle{alpha}
\bibliography{TESIS}
  
  \end{document}